\documentclass[12pt]{amsart}
\usepackage{amssymb}
\usepackage{amsfonts}
\usepackage{amsthm}
\usepackage{amsthm}
\usepackage{amscd}
\usepackage[notcite,notref]{showkeys}
\numberwithin{equation}{section}
\overfullrule=0pt \theoremstyle{plain}
\newtheorem{theorem}{Theorem}[section]

\newtheorem{proposition}[theorem]{Proposition}

\newtheorem{lemma}[theorem]{Lemma}
\theoremstyle{definition}
\newtheorem{definition}[theorem]{Definition}
\theoremstyle{remark}

\newcommand{\cM}{\mathcal{M}}

\newcommand{\Ga}{\Gamma}

\newcommand{\define}{\def}

\newcommand{\Cal}[1]{\mathcal{#1}}
\renewcommand{\frak}[1]{\mathfrak{{#1}}}
\newcommand{\refE}[1]{(\ref{E:#1})}
\newcommand{\refS}[1]{Section~\ref{S:#1}}
\newcommand{\refSS}[1]{Section~\ref{SS:#1}}

\newcommand{\refP}[1]{Proposition~\ref{P:#1}}

\newcommand{\refL}[1]{Lemma~\ref{L:#1}}

\newcommand{\C}{\ensuremath{\mathbb{C}}}
\newcommand{\N}{\ensuremath{\mathbb{N}}}

\newcommand{\Z}{\ensuremath{\mathbb{Z}}}

\renewcommand{\l}{\lambda}
\newcommand{\gb}{\overline{\mathfrak{g}}}
\newcommand{\g}{\mathfrak{g}}
\newcommand{\gh}{\widehat{\mathfrak{g}}}

\newcommand{\tr}{\mathrm{tr}}
\newcommand{\sln}{\mathfrak{sl}}
\newcommand{\gl}{\mathfrak{gl}}
\newcommand{\Hwft}{\mathcal{H}_{F,\tau}}
\newcommand{\Hwftm}{\mathcal{H}_{F,\tau}^{(m)}}

\newcommand{\KN} {Kri\-che\-ver-Novi\-kov}
\newcommand{\pN}{\ensuremath{(P_1,P_2,\ldots,P_N)}}

\newcommand{\iN}{\ensuremath{1,\ldots, N}}

\newcommand{\At}{\widetilde{\Cal A}}
\newcommand{\Lt}{\widetilde{\Cal L}}
\newcommand{\Lc}{\overline{\Cal L}}

\newcommand{\Pif} {P_{\infty}}

\newcommand{\ainf} {{\mathfrak a}_\infty} 

\newcommand{\Fln}[1][n]{F_{#1}^\lambda}

\newcommand{\A}{\mathcal{A}}

\newcommand{\V}{\mathcal{V}}

\newcommand{\Lh}{\widehat{\mathcal{L}}}
\newcommand{\eh}{\widehat{e}}

\newcommand{\kndual}[2]{\langle #1,#2\rangle}
\newcommand{\cins}{\frac 1{2\pi\mathrm{i}}\int_{C_S}}

\newcommand{\cinc}[1]{\frac 1{2\pi\mathrm{i}}\int_{#1}}

\newcommand{\w}{\omega}
\newcommand{\ord}{\operatorname{ord}}
\newcommand{\res}{\operatorname{res}}
\newcommand{\nord}[1]{:\mkern-5mu{#1}\mkern-5mu:}
\newcommand{\Fn}[1][\lambda]{\mathcal{F}^{#1}}
\newcommand{\Fl}[1][\lambda]{\mathcal{F}^{#1}}

\renewcommand{\k}{{k}}
\newcommand{\ce}{{c}}
\newcommand{\ka}{{k}}

\define\ldot{\hskip 1pt.\hskip 1pt}

\define\a{\alpha}
\redefine\d{\delta}
\define\w{\omega}

  \redefine\i{{\,\mathrm{i}}\,}
\define\ga{\gamma}
\define\cint #1{\frac 1{2\pi\i}\int_{C_{#1}}}

\define\res{\operatorname{res}}
\redefine\deg{\operatornamewithlimits{deg}}
\define\ord{\operatorname{ord}}
\define\rank{\operatorname{rank}}

\define\pfz#1{\frac {d#1}{dz}}

\define\K{\Cal K}

\redefine\H{{\mathrm{H}}}

\define\A{\Cal A}
\define\Do{\Cal D^{1}}
\define\Dh{\widehat{\mathcal{D}}^{1}}
\redefine\L{\Cal L} \redefine\D{\Cal D^{1}}
\define\KN {Kri\-che\-ver-Novi\-kov}
\define\Pif {{P_{\infty}}}
\define\Uif {{U_{\infty}}}
\define\Uifs {{U_{\infty}^*}}
\define\KM {Kac-Moody}
\define\Fln{\Cal F^\lambda_n}
\define\gb{\overline{\mathfrak{ g}}}
\define\G{\overline{\mathfrak{ g}}}
\define\Gb{\overline{\mathfrak{ g}}}
\redefine\g{\mathfrak{ g}}

\define\gh{\widehat{\mathfrak{ g}}}

\define\Lh{\widehat{\Cal L}}

\define\iN{i=1,\ldots,N}

\define\pN{p=1,\ldots,N}

\define\de{\delta}

\define\kndual#1#2{\langle #1,#2\rangle}
\define \nord #1{:\mkern-5mu{#1}\mkern-5mu:}

\define\V{\Cal V}

\begin{document}
\title[Krichever-Novikov algebras \today]
     {Affine Krichever-Novikov algebras, their representations and applications}

\thanks{Supported in part by the programme "Nonlinear dynamics and solitons" of the
Russian Academy of Science and by the Russian Foundation for Basic
Research, the project ¹ 02-01-00803 }
\author[O.K. Sheinman]{Oleg K. Sheinman}
\address[ Oleg K. Sheinman]{Steklov Mathematical Institute, ul. Gubkina, 8,
Moscow 117966, GSP-1, Russia and Independent University of Moscow,
Bolshoi Vlasievskii per. 11, Moscow 123298, Russia}
\email{oleg@sheinman.mccme.ru}

\dedicatory{Dedicated to Professor S.P.Novikov in honour of his
65th Birthday}

\begin{abstract}
The survey of the current state of the theory of Krichever-Novikov
algebras including new results on local central extensions,
invariants, representations and casimir operators.
\end{abstract}
\date{\today}
\maketitle
\tableofcontents
\section{Introduction}\label{S:intro}

In 1987, the soliton theory investigations led Krichever and
Novikov \cite{rKNFa,rKNFb,rKNFc} to the new fundamental notion in
the Lie algebra theory.  They introduced the two-dimensional
algebraic-geometrical counterpart of the celebrated Virasoro and
Kac-Moody algebras. Their definition is based on the well-known in
soliton theory type of algebraic-geometrical data, namely Riemann
surfaces with punctures and fixed (jets of) coordinates in their
neighborhoods. In case of two punctures (the only one considered
by them) this enabled them to introduce an \emph{almost graded
structure}, distinguish a unique \emph{local} central extension
among the variety of others and, finally, develop the secondary
quantization formalism (mathematically, construct modules
generated by vacuum vectors). The Krichever-Novikov algebras
\emph{contain} affine Kac-Moody and Virasoro algebras as a
subclass. This distinguish them among other (and later) approaches
to generalization of the current algebras on two dimensions.

Krichever-Novikov algebras have numerous relations to the
fundamental problems of geometry, analysis and mathematical
physics. The problem of classification of their co-adjoint orbits
turns out to be a version of the Riemann-Hilbert problem
\cite{rShea,rSha,ShSD}. In essential, these orbits are classified
by the monodromy representations of fundamental group of the
corresponding (punctured) Riemann surface. All known
representations of those algebras are parameterized (again, "in
essential") by the holomorphic vector bundles on the Riemann
surface \cite{rShf}. Thus, for the Krichever-Novikov algebras,
there exists the correspondence between orbits and representations
(originally proposed as a general principle by A.A.Kirillov). As
an essential part, this correspondence includes the well-known
correspondence between monodromy representations and holomorphic
bundles. It is interesting that, for Krichever-Novikov algebras,
the correspondence between orbits and representations has also the
form of global geometric Langlands correspondence, i.e. the
correspondence between representations of the fundamental group of
the Riemann surface of a function field (function algebra, in our
case) and of the matrix Lie algebra over this field (affine
Krichever-Novikov algebra). Below, we try to demonstrate these
relations. We also show that the the well-known Hitchin integrals
appear as co-adjoint invariants of Krichever-Novikov algebras
\cite{ShSD}, and their quantization is closely related to the
casimir-type operators which we call semi-casimirs
\cite{ShMMJ,ShUMN}.

There is a fundamental relation, based on the Kodaira-Spencer
theory, between the Krichever-Novikov algebras and the moduli
spaces of Riemann surfaces with punctures. This relation is as
follows \cite{rSSpt}: the tangent space to the moduli space of
Riemann surfaces with an arbitrary number of punctures and
arbitrary orders of fixed jets of local coordinates at punctures
is isomorphic to the direct sum of certain homogeneous subspaces
of Krichever-Novikov Virasoro-type algebra (see \cite{Kon} for the
1-puncture situation (respectively, the Virasoro algebra), and
\cite{rKNFa, GrOr} for the 2-puncture situation).

In turn, the relations between Krichever-Novikov algebras and the
moduli spaces of Riemann surfaces provide us with a basis for
certain applications in Conformal Field Theory \cite{rSSpt,
rSSpt2} which we do not concern in this publication, except for
the application to the Hitchin integrals.

The present paper is a survey of the current state of the theory
of Krichever-Novikov algebras, their invariants, representations
and some applications, including recent results.
\section{The algebras of Krichever-Novikov type}\label{S:kzkn}

Let $\Sigma$ be a compact Riemann surface of genus $g$, or in
terms of algebraic geometry, a smooth projective curve over $\C$
respectively. Let
$$
I=(P_1,\ldots,P_N),\quad  N\ge 1
$$
be a tuple of ordered, distinct points (``marked points'',
``punctures'') on $\Sigma$, $\Pif$ a distinguished point on
$\Sigma$  different from $P_i$ for every $i$, and $A=I\cup
\{\Pif\}$. The points in $I$ are called the {\it in-points}, and
the point $\Pif$ the {\it out-point}. The more general case of an
arbitrary finite set of out-points is considered in \cite{rSLc},
\cite{rSDiss}.

\subsection{The Lie algebras $\A$, $\gb$, $\L$, $\D$ and $\D_\g$}
                                   \label{SS:algs}

Let  $\A:=\A(\Sigma,I,\Pif)$ be the associative algebra of
meromorphic functions on $\Sigma$ which are regular except at the
points $P\in A$. Let  ${\g}$ be a complex finite-dimensional
reductive Lie algebra. Then
\begin{equation}\label{E:curalg}
   {\gb}=\g\otimes_{\mathbb C}{\A}
\end{equation}
is called the {\it Krichever-Novikov current algebra}
\cite{rKNFa,rShea}. The Lie bracket on $\gb$ is given by the
relations
\begin{equation}\label{E:curr}
[x\otimes A,y\otimes B]=[x,y]\otimes AB.
\end{equation}
We will often suppress the symbol $\otimes$ in our notation.

Let $\L$ denote the Lie algebra  of meromorphic vector fields on
$\Sigma$ which are allowed to have poles only at the points $P\in
A$ \cite{rKNFa,rKNFb,rKNFc}.

For the Riemann sphere ($g=0$) with quasi-global coordinate $z$,
$I=\{0\}$ and $\Pif=\infty$, $\A$ is the algebra of Laurent
polynomials, the  current algebra $\gb$ is the loop algebra, and
$\L$ is the Witt algebra.

The algebra $\L$  operates on the elements of $\A$ by the (Lie)
derivative. This allows to define the Lie algebra of differential
operators $\Do$ of degree $\le 1$ as their  semi-direct {sum}. As
a vector space, $\Do=\A\oplus \L$. The Lie structure is defined by
\begin{equation}
[(g,e),(h,f)]:=(e\ldot h-f\ldot g,[e,f]),\quad g,h\in\A,\
e,f\in\L.
\end{equation}
Here $e\ldot f$ denotes  the (Lie) derivative.

Similarly, we define the Lie algebra  $\D_\g$. As a vector space,
$\D_\g=\gb\oplus\L$. The Lie structure on $\gb$, $\L$ is as above,
and additionally
\begin{equation}\label{E:comm}
  [e,x\otimes A]:=-[x\otimes A,e]:=x\otimes (e\ldot A).
\end{equation}
In particular, for $\g={\frak{gl}}(1)$ one obtains $\D_\g=\D$ as a
special case.
\subsection{Meromorphic forms of weight $\lambda$ and
\KN\ duality}\label{SS:mforms}
Let $\K$ be the canonical line bundle. For every $\l\in\Z$ we
consider the bundle $\ \K^\l:=\K^{\otimes \l}$. Here we use the
usual convention: $\K^0=\Cal O$ is the trivial bundle,  and
$\K^{-1}=\K^*$ is the holomorphic tangent line bundle. Indeed,
after fixing a theta characteristics, i.e. a bundle  $S$ with
$S^{\otimes 2}=\K$, it is possible to consider $\l\in \frac
{1}{2}\Z$. Denote by $\Fl$ the (infinite-dimensional) vector space
of global meromorphic sections of $\K^\l$ which are holomorphic on
$\Sigma\setminus A$. The elements of $\Fl$  are called
(meromorphic) forms or tensors of weight $\lambda$.

The cases of a special interest are as follows: functions
($\l=0$), vector fields ($\l=-1$), 1-forms ($\l=1$), and
quadratic differentials ($\l=2$).  The space of functions is
already denoted by $\A$, and the space of vector fields by $\L$.

Denote by $C_S$ any cycle homologous to a small circle surrounding
$P_\infty$. Following \cite{rKNFa}, we refer to any $C_S$ as to a
\emph{separating cycle}. The notion of separating cycle is
introduced by Krichever and Novikov and closely related to their
conceptions of locality and almost-grading.
\begin{definition}\label{D:knpair}
The {\it Krichever-Novikov pairing} ({\it KN pairing}) is the
pairing between $\Fl$ and $\Fl[1-\l]$ given by
\begin{equation}\label{E:knpair}
\begin{gathered}
\Fl\times\Fn[1-\l]\ \to\ \C,
\\
\kndual {f}{g}:=\cins f\otimes g =\sum_{P\in I}\res_{P}(f\otimes
g)= -\res_{\Pif}(f\otimes g),
\end{gathered}
\end{equation}
where $C_S$ is any separating cycle.
\end{definition}
The last equality follows from the residue theorem. Observe that
in \refE{knpair} the integral does not depend on the separating
cycle chosen.

\subsection{\KN\ bases}\label{SS:knb}

In \cite{rKNFa}, Krichever and Novikov introduced special bases
for the spaces of meromorphic  tensors on Riemann surfaces with
two marked points. For $g=0$ the Krichever-Novikov bases coincide
with the Laurent bases. The multi-point generalization of these
bases is given in \cite{rSLc,rSDiss} (see also \cite{rSad}). We
define here the Krichever-Novikov type bases for the tensors of an
arbitrary weight $\lambda$ on Riemann surfaces with $N$ marked
points as introduced in \cite{rSLc,rSDiss}.

For fixed $\l$ and for  every $n\in\Z$, and $p=1,\ldots,N$ we
exhibit a certain element $f_{n,p}^\l\in\Fl$. The basis elements
are chosen in such a way that they fulfill the duality relation
\begin{equation}\label{E:edu}
\kndual {f_{n,p}^\l} {f_{m,r}^{1-\l}}= \de_{-n}^{m}\cdot
\de_{p}^{r}
\end{equation}
with respect to the KN pairing \refE{knpair}. In particular, this
implies that the KN pairing is non-degenerate. Additionally, the
elements fulfill
\begin{equation}\label{E:ordfn}
\ord_{P_i}(f_{n,p}^\l)=(n+1-\l)-\d_i^p,\quad i=1,\ldots,N .
\end{equation}
The recipe for choosing the order at the point $\Pif$ is such that
up to a scalar multiplication there is a unique such element which
also fulfills \refE{edu}. After choosing local coordinates $z_p$
at the points $P_p$ the scalar can be fixed by requiring
\begin{equation}
{f_{n,p}^\l}(z_p)=z_p^{n-\l}(1+O(z_p))\left(dz_p\right)^\l, \quad
p=1,\ldots,N.
\end{equation}
To give an impression about requirement at $P_\infty$ let us
consider the case $g\ge 2$, $\l\ne 0,1$ and a generic choice for
the points in $A$ (or $g=0$ without any restriction). Then we
require
\begin{equation}\label{E:ordi}
\ord_{\Pif}(f_{n,p}^\l)=-N\cdot(n+1-\l) +(2\l-1)(g-1)\ .
\end{equation}
By Riemann-Roch type arguments, it is shown in \cite{rSLa} that
there exists only one such element.

Explicit expressions for the basis elements $f_{n,p}^\l$ in terms
of rational functions ($g=0$), the Weierstra\ss\ functions
($g=1$), prime forms and theta functions ($g\ge1$) are given in
\cite{rKNFa}, \cite{rSLb}, \cite{rRDS}, \cite{rSDeg}. For $g=0$
and $g=1$, such expressions can be found also in \cite[\S\S
2,7]{rSSpt}. One can see from the explicit expressions that the
basis elements ``analytically'' depend on the complex structure of
the Riemann surface.

For the following cases we introduce a special notation:
\begin{equation}\label{E:econc}
A_{n,p}:=f_{n,p}^0,\quad e_{n,p}:=f_{n,p}^{-1},\quad
\w^{n,p}:=f_{-n,p}^1,\quad \Omega^{n,p}:=f_{-n,p}^2 \ .
\end{equation}

For $g=0$ and $N=1$ the basis elements \refE{econc} coincide with
the standard generators of the Laurent polynomials, the Witt
algebra and their dual spaces, respectively. For $g\ge 1$ and
$N=1$ these elements coincide, up to a shift of index, with those
given by Krichever and Novikov \cite{rKNFa,rKNFb,rKNFc}.

\subsection{Almost graded structure, triangular
  decompositions}\label{SS:ags}${ }$
For $g=0$ and $N=1$ the Lie algebras introduced in \refSS{algs}
are graded. A grading is a necessary tool for developing their
structure theory and the theory of their highest weight
representations. For the higher genus case (and for the
multi-point situation for $g=0$) there is no grading. It is a
fundamental observation due to Krichever and Novikov
\cite{rKNFa,rKNFb, rKNFc} that a weaker concept, an almost
grading, is sufficient to develop a suitable structure and
representation theory in this more general context.

An (associative or Lie) algebra is called {\it almost-graded} if
it admits a direct decomposition as a vector space $\
\V=\bigoplus_{n\in\Z} \V_n\ $, where (1) $\ \dim \V_n<\infty\ $
and (2) there are constants $R$ and  $S$ such that
\begin{equation}\label{E:eaga}
 \V_n\cdot \V_m\quad \subseteq \bigoplus_{h=n+m-R}^{n+m+S} \V_h,
 \qquad\forall n,m\in\Z\ .
\end{equation}
The elements of $\V_n$ are called {\it homogeneous  elements of
degree $n$}. Let $\V=\oplus_{n\in\Z} \V_n$ be an almost-graded
algebra and
 $M$ an $\V$-module.
The module $M$ is called an {\it almost-graded} $\V$-module if it
admits a direct decomposition as a vector space $\
M=\bigoplus_{m\in\Z} M_m\ $, where (1) $\ \dim M_m<\infty\ $ and
(2) there are constants $T$ and  $U$ such that
\begin{equation}\label{E:eagm}
 V_n\ldot M_m\quad \subseteq \bigoplus_{h=n+m-T}^{n+m+U} M_h,
 \qquad\forall n,m\in\Z\ .
\end{equation}
The elements of $M_n$ are called {\it homogeneous  elements of
degree $n$}.

For the space $\Fl$ of \refSS{mforms}, its homogeneous degree $n$
subspace $\Fln$ is defined as the subspace generated by the
elements $f_{n,p}^\l$, $p=1,\ldots,N$. Then
$\Fl=\bigoplus_{n\in\Z}\Fln$.
\begin{proposition}\label{P:almgrad}\cite{rSLc,rSDiss}
With respect to the introduced degree, the vector field algebra
$\L$, the function algebra $\A$, and the differential operator
algebra $\D$  are almost-graded and the $\Fl$ are  almost-graded
modules over them.
\end{proposition}
\noindent

For particular algebras, the almost graded structure has a more
special description.
\begin{proposition}\cite{rSLc,rSDiss} \label{P:struct}
There exist constants $\ K,L,M\in\N\ $ such that for all
$n,m\in\Z$
\begin{equation}\label{E:struct}
\begin{aligned}
A_{n,p}\cdot A_{m,r}&=\d_p^r\,A_{n+m,p}+
\sum_{h=n+m+1}^{n+m+K}\sum_{s=1}^N\alpha_{(n,p),(m,r)}^{(h,s)}
A_{h,s}\ ,
\\
[e_{n,p},e_{m,r}]&=\d_p^r\,(m-n)\,e_{n+m,p}+
\sum_{h=n+m+1}^{n+m+L}\sum_{s=1}^N\gamma_{(n,p),(m,r)}^{(h,s)}
e_{h,s}\ ,
\\
e_{n,p}\ldot A_{m,r}&=\d_p^r m \,A_{n+m,p}+
\sum_{h=n+m+1}^{n+m+M}\sum_{s=1}^N\beta_{(n,p),(m,r)}^{(h,s)}
A_{h,s}\ ,
\end{aligned}
\end{equation}
with suitable coefficients $\ \alpha_{(n,p),(m,r)}^{(h,s)},
\beta_{(n,p),(m,r)}^{(h,s)}, \gamma_{(n,p),(m,r)}^{(h,s)}\in\C$.
\end{proposition}
The constants $K,L$ and $M$ depend on the genus $g$ and the number
of points $N$. Explicit expressions  can be found in \cite[\S
2]{rSSpt}.

As a vector space, the algebra  $\A$ can be decomposed as follows:
\begin{equation}\label{E:edeca}
\begin{gathered}
\A=\A_+\oplus\A_{(0)}\oplus\A_-,
\\
\A_+:=\langle A_{n,p}\mid n\ge 1,\pN\rangle,\quad \A_-:=\langle
A_{n,p}\mid n\le -K-1,\pN\rangle\ ,\quad
\\
\A_{(0)}:=\langle A_{n,p}\mid -K\le  n\le 0, \pN\rangle \ .
\end{gathered}
\end{equation}
and the Lie algebra $\L$ as follows:
\begin{equation}\label{E:edecv}
\begin{gathered}
\L=\L_+\oplus\L_{(0)}\oplus\L_-,
\\
\L_+:=\langle e_{n,p}\mid n\ge 1,\pN\rangle,\quad \L_-:=\langle
e_{n,p}\mid n\le -L-1,\pN\rangle\ ,\quad
\\
\L_{(0)}:=\langle e_{n,p}\mid -L\le  n\le 0, \pN\rangle \ .
\end{gathered}
\end{equation}
We call \refE{edeca}, \refE{edecv} the {\it triangular
decompositions}. In a similar way we obtain a triangular
decomposition of $\D$.

Due to the almost-grading the subspaces $\A_\pm$ and $\L_\pm$ are
subalgebras but the subspaces $\A_{(0)}$, and $\L_{(0)}$, in
general, are not. We use the term {\it critical strip} for  the
latters.

Observe that $\A_+$ and $\L_+$ can be described as the algebra of
functions (respectively, vector fields) having a zero of, at
least, order one (respectively, two) at the points $P_i,\iN$.
These algebras can be enlarged by adding all elements which are
regular at all $P_i$'s (these are $\ \{A_{0,p},\pN\}\ $,
respectively $\ \{e_{0,p},e_{-1,p},\iN\}\ $). We denote the
enlarged algebras by $\A_+^*$, resp. by $\L_+^*$.

On the other hand $\A_-$ and $\L_-$ could also be enlarged such
that they contain all elements  which are regular at $\Pif$. This
is explained in detail in \cite{rSSpt}. We obtain  $\A_-^*$ and
$\L_-^*$ respectively. For every $p\in\N_0$, let $\L_-^{(p)}$ be
the subalgebra of vector fields of order $\ge p+1$ at the point
$\Pif$, and $\A_-^{(p)}$ be the subalgebra of functions of order
$\ge p$ at the point $\Pif$. We obtain a decomposition
\begin{equation}
\L=\L_+\oplus\L_{(0)}^{(p)}\oplus\L_-^{(p)},\ \text{ for $p\ge 0$}
\quad \text{and}\quad
\A=\A_+\oplus\A_{(0)}^{(p)}\oplus\A_-^{(p)},\ \text{ for $p\ge
1$},
\end{equation}
with ``critical strips'' $\L_{(0)}^{(p)}$ and $\A_{(0)}^{(p)}$,
which are only subspaces. The case of $p=1$ is of particular
interest. We call $\L_{(0)}^{(1)}$ the \emph{reduced critical
strip}. For $g\ge 2$ its dimension is equal to
\begin{equation}\label{E:dimcs}
\dim \L_{(0)}^{(1)}=N+N+(3g-3)+1+1=2N+3g-1\ .
\end{equation}
On the right hand side, the first two terms correspond to $\L_0$
and $\L_{-1}$. The intermediate term comes from the vector fields
in the basis which have poles at the $P_i,\iN$ and $\Pif$. The
$1+1$ corresponds to the basis vector fields with  exact order
zero (one) at $\Pif$.

On the {\it higher genus current algebra} $\gb$, the
almost-grading is introduced by setting  $\deg(x\otimes
A_{n,p}):=n$. As above, we obtain a triangular decomposition
\begin{equation}\label{E:espcur}
\G=\G_+\oplus\G_{(0)}\oplus\G_-,\quad\text{with}\quad
\G_\beta=\g\otimes \A_\beta,\quad \beta\in\{-,(0),+\}\ ,
\end{equation}
In particular, $\G_\pm$ are subalgebras. The corresponding is true
for the enlarged subalgebras. Among them, $
\gb_{reg}:=\gb_-^{(1)}=\g\otimes \A_-^{(1)} $ is of special
importance. It is called the regular subalgebra.

The finite-dimensional Lie algebra $\g$ can be naturally
considered as a subalgebra of $\Gb$. It is a subspace of $\G_0$.
To see this, make use of the relation $1=\sum_{p=1}^N A_{0,p}$
\cite[Lemma 2.6]{rSSpt}.


\subsection{Local central extensions}\label{SS:cocyc}

Let $\V$ be a Lie algebra and $\ga$  a 2-cocycle on $\V$, i.e. an
antisymmetric bilinear form obeying
\begin{equation}\label{E:cocycle}
\ga([f,g],h)+\ga([g,h],f)+\ga([h,f],g)=0,\quad\forall f,g,h\in\V.
\end{equation}
Given $\ga$, we can define a Lie algebra structure on
$\widehat\V=\V\oplus\C t$ by
\begin{equation}
[f,g]\,\hat{}:=[f,g]+\ga(f,g)\cdot t,\quad
[t,\widehat\V]\,\hat{}=0
\end{equation}
where $t$ is a formal central generator. Up to equivalence,
central extensions are classified by the elements of
$\H^2(\V,\C)$, the second Lie algebra cohomology space with values
in the trivial module $\C$. In particular, two cocycles
$\ga_1,\ga_2$ define equivalent central extensions if and only if
there exist a linear form $\phi$ on $\V$ such that
\begin{equation}
\ga_1(f,g)=\ga_2(f,g)+\phi([f,g]).
\end{equation}

\begin{definition}\label{D:local}
Let $\V=\bigoplus_{n\in\Z} \V_n$ be an almost-graded Lie algebra.
A cocycle $\ga$ for $\V$ is called local (with respect to the
almost-grading) if there exist $M_1,M_2\in\Z$ with
\begin{equation}\label{E:local}
\forall n,m\in\Z:\quad \ga(\V_n,\V_m)\ne 0\implies M_2\le n+m\le
M_1.
\end{equation}
\end{definition}
Assuming $\deg(t)=0$, \emph{the central extension $\widehat\V$ is
almost-graded if and only if it is given by a local cocycle
$\ga$}. In this case we call $\widehat\V$ an {\it almost-graded
central extension} or a {\it local central extension}.

For $g=0$ the description of 2-cocycles is due to
I.Gelfand-D.Fuchs, V.Kac, E.Arbarello-C.DeConcini-V.Kac-C.Procesi.
For $g>0$ the problem admits a direct generalization, namely, the
problem of description of the \emph{universal} central extensions.
For different Lie algebras on Riemann surfaces such description is
given by I.Frenkel, E.Frenkel, P.Etingof, Ch.Kassel, M.Bremner,
D.Millionschikov. As it is mentioned above, only local central
extensions enable one to develop the highest weight (physically,
the secondary quantization) formalism. The problem of
classification of the \emph{local} central extensions reduces to
the description of local 2-cocycles. This problem originates in
the work of I.Krichever and S.Novikov \cite{rKNFa} where the
authors introduced what we call below the \emph{geometric
cocycles} for $\A$ and $\L$ and, for $\L$, outlined the proof of
the 1-dimensionality of the space of local cocycles.

The full classification of the local cocycles for $\gb$ and
$\Do_\g$ with a reductive $\g$, including full proofs, is given by
M.Schlichenmaier \cite{SchlCo, SchlMMJ}. For the Lie algebras of
\refS{kzkn}, he has proven that \emph{any local 2-cocycle is
cohomologous to a linear combination of the basic geometric
cocycles}. Below, we give the full list of them. The more precise
result of \cite{SchlCo, SchlMMJ} for the Lie algebra $\gb$ is as
follows: \emph{any local 2-cocycle is cohomologous to an
$\L$-invariant one, and the latter is geometric} (by definition, a
cocycle $\ga$ on $\gb$ is \emph{$\L$-invariant} if
$\ga(x(e.A),yB)=\ga(x(A),y(e.B))$ for any $e\in\L$, $x,y\in\g$,
$A,B\in\A$).

The list of basic geometric cocycles is as follows:

--- for $\A$ \cite{rKNFa}:
\begin{equation}\label{E:fung}
\ga^{(\A)}(g,h):=\cinc{C_S} gdh ;
\end{equation}

--- for $\L$ \cite{rKNFa}:
\begin{equation}\label{E:ecva}
\ga^{(\L)}(e,f) :=\cinc{C_S}\left(\frac 12(\tilde e'''\tilde
f-\tilde e\tilde f''') -R\cdot(\tilde e'\tilde f-\tilde e\tilde
f')\right)dz\
\end{equation}
where $e=\tilde e\frac {d}{dz}$ and $f=\tilde f\frac {d}{dz}$ are
the local representations of the involved vector fields, $R$ be a
global holomorphic projective connection (see \cite{rKNFa,rSSpt}
for details). A different choice of the projective connection
(even if we allow meromorphic projective connections with poles
only at the points in $A$) yields a cohomologous cocycle, hence an
equivalent central extension;

--- for $\Do$: \refE{fung} \refE{ecva} and the following \emph{mixing cocycle}
\begin{equation}\label{E:mixg}
\ga^{(m)}(e,g):=-\ga^{(m)}(g,e):= \cinc {C_S}\left(\tilde e\cdot
g''+T\cdot (\tilde e\cdot
 g')\right)dz
\end{equation}
where $T$ is an \emph{affine connection} (see detailes in
\cite{rSDiss}, \cite{ShMMJ}, \cite{SchlCo}) which is holomorphic
outside $A$ and has at most a pole of order one at $\Pif$. Again,
the cohomology class does not depend on the chosen affine
connection;

--- for $\gb$:
\begin{equation}\label{E:coaff}
\ga^{(\g)}(x\otimes f,y\otimes g)=\a(x,y)\ga^{(\A)}(g,h)
\end{equation}
where $x,y\in\g$, $f,g\in\A$, $\a$ is a bilinear form. We give
here the full list of basic geometric cocycles only for
$\g={\sln}(n)$ and $\g=\gl(n)$. In the first case, it consists
from only one cocycle of the form \refE{coaff} corresponding to
the (Cartan-Killing) form $\a(x,y)=\tr(ad\, x\cdot ad\, y)$. In
the second case, there are two independent cocycles corresponding
to $\a_1(x,y)=\tr(xy)$ and $\a_2(x,y)=\tr(x)\tr(y)$, respectively;

--- for $\Do_\g$, $\g=\gl(n)$:
\[\ga^{(m,\g)}(e,x\otimes g)=\tr(x)\ga^{(m)}(e,g)\]
where again $x\in\g$, $e\in\L$, $g\in\A$. Hence, $\ga^{(m,\g)}$
for $\g=\sln(n)$. The absence of mixing cocycles in semi-simple
case was observed already in \cite{ShMMJ}.

For the estimation of boundaries of locality and other details we
refer to \cite{SchlCo, SchlMMJ}.


\subsection{Affine algebras}\label{SS:affine}

Now we are in position to define affine Krichever-Novikov
algebras. Let $\g$ be a reductive finite-dimensional Lie algebra.
By \emph{affine Krichever-Novikov algebra} we mean the central
extension of $\gb$ given by a cocycle of the form \refE{coaff}.

Thus, as a linear space, an affine Krichever-Novikov algebra is
$\gh=\g\otimes\A\oplus\C t$, and the Lie bracket is given by
\begin{equation}\label{E:eaff}
[x\otimes f,y\otimes g]= [x,y]\otimes (f
g)+\alpha(x,y)\cdot\gamma^{(\A)}(f,g)\cdot t,\qquad [\,t,\gh]=0\ ,
\end{equation}
where $\gamma^{(\A)}$ is given by \refE{fung}, $x,y\in\g$,
$f,g\in\A$, $t$ is a formal central generator. Certainly, $\gh$
depends on the bilinear form $\alpha$. For $\g={\sln}(n)$ and
$\g=\gl(n)$ we take $\a(x,y)=\tr(xy)$.

The cocycle defining the central extension is local, hence $\gh$
is almost-graded (we set $\deg t:=0$ and $\deg(\widehat{x\otimes
A_{n,p}}):=n$). Again we obtain a triangular decomposition
\begin{equation}\label{E:espaff}
\gh=\gh_+\oplus\gh_{(0)}\oplus\gh_-\quad\text{with}\quad
\gh_\pm\cong\g_\pm\quad\text{and}\quad\gh_{(0)}=\g_{(0)}\oplus\C\cdot
t\ .
\end{equation}
The corresponding is true for the enlarged algebras. Among them,
$$
\gh_{reg}:=\gh_-^{(1)}=\gb_-^{(1)}=\g\otimes \A_-^{(1)},\qquad
\gh_+^{*,ext}=\gb_+^*\oplus \C \;t=(\g\otimes \A_+^*)\oplus\C \;t\
.
$$
are of special interest.


\section{Orbits and invariants}\label{S:orb}

In this section, following the lines of \cite{rShea,rSha,ShSD}, we
give a description of the co-adjoint orbits and invariants of an
affine Krichever-Novikov algebra ${\gh}$.  For simplicity, we
consider the case of the Riemann surface with two marked points.
Thus, ${\gh} =\g\otimes_{\C}{\A}\oplus{\C}t $ where
$\A=\A(\Sigma,P_\pm)$, $\Sigma$ is a compact algebraic curve over
$\C$ with two marked points $P_{\pm}$ ($genus\,\Sigma
>0$), $\A(\Sigma,P_\pm)$ stays for the algebra of meromorphic
functions on $\Sigma$ which are regular out the points $P_\pm$,
${\g}$ is a complex reductive Lie algebra.

Let $\omega$ be a meromorphic ${\g}^*$-valued 1-form on $\Sigma$
regular outside the points $P_\pm$, $b\in{\Bbb C}$. We identify
the dual space ${\gh}^*$ with the space of operators of the form
$bd+\omega$ which take the sections of the $l$-dimensional trivial
bundle on $\Sigma$ to the sections of the tensor product of the
same bundle by $\K$. Here, $d$ denotes the differential on
$\Sigma$.

Consider the mappings $\Sigma\to\exp\g$ holomorphic outside
$P_\pm$. Such mappings, with the operation of a point-wise
multiplication, form the group which we denote $G$. We call $G$
the \emph{current group}, and its elements the \emph{group
currents}. The current group acts on the elements of ${\gh}^*$ by
{\it gauge transformations}: $\omega\mapsto g\omega g^{-1}-b\cdot
dg\cdot g^{-1}$, $b\mapsto b$. We want to describe the space of
orbits of this action. Since we made no assumption on the behavior
of the group currents at $P_\pm$, in general, ${\gh}^*$ is not
invariant with respect to the gauge action. There are two ways
out. The first is to consider orbits in some bigger space. It is
an interesting question what space it should be. Our
\emph{conjecture} is that this is a space of matrix
Backer-Achiezer functions. Here, we choose another approach. By
\emph{orbit} we mean the intersection of the true orbit with
${\gh}^*$.
\begin{theorem}\label{T:orbs}\cite{rShea,rSha}
The space of orbits in generic position which are in an invariant
affine hyperplane $b=const$ ($b\ne 0$) in ${\gh}^*$ is in a
one-to-one correspondence with the space of equivalence classes of
representations $\tau :\pi_1(\Sigma\backslash P_\pm
)\rightarrow\exp\g$ such that $\tau(\gamma)$ is a semisimple
element (where $\gamma$ is a homotopy class of the separating
contour $C_S$, see \refSS{knb}).
\end{theorem}

We call $\tau$ the monodromy representation of
$\pi_1(\Sigma\backslash P_\pm )$. The correspondence between the
orbits and the monodromy representations is established as follows
\cite{rShea, rSha}: each element $bd+\omega\in\gh^*$ is assigned
with the \emph{monodromy equation}
\begin{equation}\label{E:monoe}
(bd+\omega)g=0,
\end{equation}
and, further on, with the monodromy representation of this
equation. For gauge equivalent elements, this procedure results in
the equivalent monodromy representations and vise verse.

To complete the classification of the orbits (at least, in a
generic position), we should check that \emph{any} monodromy
representation (of a certain class) corresponds to some equation
of the form \refE{monoe}, i.e. resolve the \emph{Riemann-Hilbert
problem}. Under assumption of semi-simplicity of $\tau(\gamma)$,
this problem can be resolved which follows from the results of
\cite{Arn}, see also \cite[Proposition 1.5.2]{Bol}.

Define \emph{co-adjoint invariants} of $\gh$ as such functions on
$\gh^*$ that are constant on the co-adjoint orbits, hence can be
pushed down to the orbit space. It is our next step to construct a
full system of independent co-adjoint invariants. We will do it
for $\g=\gl(n)$.

Following \cite{Bol}, each monodromy representation coupled with a
set of highest weights of $\g$ assigned to the marked points can
be associated with a smooth vector bundle on $\Sigma$ and a flat
logarithmic (see below for the definition) connection $\nabla$ on
this bundle. We interpret this as an assignment of a bundle (and a
connection) to each co-adjoint \emph{orbit} of $\gh$. It is our
conjecture that the additional parameters (highest weights)
correspond to the above mentioned hypothetical extension of the
dual space $\g^*$.

Going over to the construction of the invariants, take an
arbitrary dominant weight $\l_+=(\l_1,\ldots,\l_n)$
($\l_1\ge\l_2\ge\ldots\ge\l_n$) at the point $P_+$ and a dominant
weight $\l_-$ of the form $(*,0,\ldots,0)$ or $(0,\ldots,0,*)$ at
the point $P_-$. Denote by $B_{\l_+,\l_-,\tau}$ the corresponding
vector bundle on $\Sigma$. In accordance with \cite{Bol}, the
weight $\l_-$ of such form is uniquely determined by the triple
consisting of $\l_+$, $\tau$ and the degree of the bundle. Let us
require the degree to be equal to $gn$, ($g=genus(\Sigma)$,
$n=\rank\g$), then we can suppress $\l_-$ in the notation of the
bundle and write it down as $B_{\l,\tau}$ where $\l=\l_+$.

According to the Donaldson theorem \cite{Don,Had}, $\nabla$ is
gauge equivalent to the connection of the form $d+A+\phi+\phi^*$
where $(A,\phi)$ is a solution to the self-duality equations
$d''_A\phi=0$, $F(A)+[\phi,\phi^*]=0$, $F(A)$ is the curvature of
the connection $A$, $*$ is the Hermitian conjugation (for $\g=\g
l(n)$ ), $d+A=d'_A+d''_A$ is the standard notation for a
connection in complex coordinates, and the form of the connection
$A$ has no singularities at the marked points. Define the complex
structure on $B_{\l,\tau}$ by means of the $d''_A$ as a
$\overline\partial$-operator. According to \cite{Had}, the bundle
$B_{\l,\tau}$ with this complex structure is a holomorphic bundle
on $\Sigma$ (in particular this implies
$\phi\in\Omega^{1,0}(\Sigma\backslash P_\pm)$,
$\phi^*\in\Omega^{0,1}(\Sigma\backslash P_\pm)$, hence, both these
forms can be retrieved from the canonical form of the connection).
The just constructed holomorphic bundle $B_{\l,\tau}$ has the rank
$n$ and the degree $gn$. We assume that $B_{\l,\tau}$ is generic
in sense of \cite{rKNU} (see also \cite{rShf}), hence {\it the
space of its meromorphic sections which are holomorphic out the
marked points $P_\pm$ is isomorphic to the space $\Cal F$ of
meromorphic vector-functions $f=(f_1,\ldots,f_l)^T$ on the same
Riemann surface which are holomorphic out the points $P_\pm$ and
the divisor $D=\Sigma_{i=1}^{gn}\gamma_i$ of degree $gn$, have at
most simple poles at the points of the latter and satisfy the
following relations
   $$ (\res_{\ga_i}f_j)\a_{ik}=(\res_{\ga_{i}}f_k)\a_{ij},\,\,\,
      i=1,\ldots,gl,\,\,\, j=1,\ldots,l,
   $$
where $\a_{ik}$ are constants}. The points of the divisor $D$ and
the numbers $\a_{ik}$ are called {\it the Tjurin parameters} of
the bundle.

Our goal is to find the invariants of the gauge action on $\gh^*$.
To this end, first, take $\l=0$. A self-dual pair $A,\phi$ has the
well-known invariants called \emph{Hitchin integrals} \cite{rHsb},
namely the coefficients of the expansion of $\tr\,\phi^k$ over the
basis holomorphic tensors of weight $k$, $k=1,\ldots, n$. For
$\g=\gl(n)$, it follows from \cite{rHsb} that their number is
equal to $n^2(g-1)+1$. Provided the value of $\tau$ on the class
of the separating contour is nontrivial, $\phi$, generically, has
simple poles at the points $P_\pm$. The space of such forms has
the  dimension $n^2(g+1)$ (the contribution of holomorphic forms
is equal to $n^2g$, and $n^2$ is contributed by residues of $\phi$
at $P_\pm$, taking into account the relation $\res_{P_+}\phi
+\res_{P_-}\phi=0$), hence, any such form is determined by
$n^2(g+1)$ invariants.

To obtain these invariants, the basis holomorphic tensors should
be supplemented with tensors that have poles up to the order $n^2$
at $P_\pm$. The number of independent main parts of such tensors
is $2n^2$ ($n$ weights, $n$ orders of poles and two points). Thus,
we have $2n^2$ singular basis elements in addition to regular
ones, which gives, a'priori, $n^2(g-1)+1+2n^2=n^2(g+1)+1$
integrals in total. The tensor $\tr\,\phi^k$ has non-trivial
projections on (singular) basic tensors of weight $k$ and order
$\le k$. This gives $2k$ integrals for each $k$. To obtain some
non-trivial projections on the basis tensors of weight $k$ and
orders $k+1,\ldots,n$, we take the expansions of $\tr\,A_m\phi^k$,
$m=-g-n+k,\ldots,-g-1,1,\ldots,n-k$. These $A_m$'s have orders
$1,\ldots,n-k$ at $P_\pm$. They span the subspace of dimension
$2(n-k)$ which gives $2(n-k)$ additional (singular) integrals. In
total, we have $2n$ singular integrals for any $k=1,\ldots,n$,
thus, $2n^2$ integrals in addition to the regular (Hitchin)
integrals. Observe that in the original Hitchin situation (only
holomorphic tensors are considered) we should introduce only
global holomorphic functions (i.e. constants), thus obtain no new
integrals. By the same reason, for $k=n$ all integrals come from
$\tr\,\phi^n$. Since $\res_{P_+}\phi +\res_{P_-}\phi=0$, the
coefficients in the terms of degree $n$ of the Laurent expansions
of the tensor $\tr\,\phi^n$ at the marked points either are equal
or differ by a sign (depending on parity of $n$). By this, one
parameter is eliminated. Thus, we obtain $n^2(g+1)$ {\it
generalized Hitchin invariants}.

One more set of invariants, which depend (via the holomorphic
structure) on the component $A$ of the self-dual pair, is the
tuple of above introduced Tjurin parameters of the bundle
$B_{0,\tau}$ \cite{rTur,rKNU}. There is $n^2(g-1)+1$ independent
ones among them. In addition, there is one "trivial" invariant $b$
(see above). In total, we have constructed $2n^2g+2$ independent
invariants. As it is easy to show, this number equals to the
dimension of the orbit space. Thus {\it the generalized Hitchin
invariants and the Tjurin invariants form the the full family of
independent invariants of the gauge action on $\gh^*$}.


\section{Representations of \KN\ algebras}\label{S:kzrep}
Let $B$ be a holomorphic vector bundle on $\Sigma$, $\nabla$ a
meromorphic connection on $B$, $\tau$ an irreducible
representation of $\g$. For this section, our goal is to assign
each triple $(B,\nabla,\tau)$ with a wedge representation of
$\Do_\g$ and investigate the properties of this representation as
a $\gh$-module.

\subsection{Krichever--Novikov bases for the holomorphic vector
bundles}\label{SS:knbb}

Actually, these are the bases in the function space $\mathcal F$
(see \refS{orb}) which is isomorphic to the space of meromorphic
sections of $B$. The bases in question where introduced in
\cite{rKNR2p} (for $N=1$) and then applied to the construction of
wedge (or \emph{fermion}) representations in \cite{rShf}. Here,
combining the approaches of \cite{rKNR2p} and \refSS{knb} (going
back to \cite{rSLa,rSLb}) we introduce these bases for an
arbitrary $N$.

From now on, let $n$ stay for the degree of the Krichever-Novikov
basis element (to be defined below), $l$ stay for the rank of $\g$
and $r$ stay for the rank of $B$.

For any tryple of integers $n\in\Z$, $j=0,1,\ldots r-1$ and
$p=1,\ldots, N$, let $\psi_{n,j,p}=(\psi_{n,j,p}^i)$ be a vector
function ($i=0,1,\ldots, r-1$) in $\mathcal F$. This function is
specified by its asymptotic behavior at the points of $A$ which is
assumed to be as follows:
\begin{equation}\label{E:bas1}
  \psi_{n,j,p}^i(z_q)=z_q^{n+1-\d_{p,q}}(\xi_{npqj}^i+O(z_q))
\end{equation}
where $z_q$ is a local coordinate at $P_q$, $q\in\{1,\ldots,N\}$,
$\d_{p,q}$ is the Kronecker symbol, $\xi_{npqj}$ are the complex
numbers such that $\xi_{nppi}^i=1$,\ \ $\xi_{nppj}^i=0$, $i>j$;
\begin{equation}\label{E:bas2}
  \psi_{n,j,p}^i(z_\infty)=z_\infty^{-nN-N+1}(\xi_{np\infty j}^i+O(z_\infty))
\end{equation}
where $\xi_{np\infty j}^i=0$, $i<j$.

For a given $n$, the space spanned by vector-valued functions
$\psi_{n,j,p}$ is referred to as the space of degree $n$
functions.

\begin{proposition}\label{P:2.2}
1. There exists a unique vector-valued function $\psi_{n,j,p}$
satisfying conditions \refE{bas1} and \refE{bas2}.

2. For a given $n$, the dimension of the space spanned by the
vector-valued functions $\psi_{n,j,p}$ is equal to~$rN$.
\end{proposition}

\begin{proof}
Introduce the divisor $D_{n,p}=D+nP_p+\Sigma_{q\ne p}
(n+1)P_q+(-nN-N+1)P_\infty$, hence $\deg\,D_{n,p}=\deg\,D=-rg$.
Let $(D_{n,p})$ denote the space of (scalar-valued) meromorphic
functions with a divisor not less than  $D_{n,p}$, then $\dim
(D_{n,p})=rg-g+1$. For the analogous space of functions taking
values in an $r$-dimensional vector space (denote it by
$(D_{n,p})_r$), we have $\dim\, (D_{n,p})_r=r(rg-g+1)$. Such
functions in $\mathcal F$ satisfy the $(r-1)gr$ Tjurin relations.
Therefore, $\dim\, ((D_{n,p})_r\cap {\mathcal F}) =
r(gr-g+1)-(r-1)gr=r$.

Observe that $\psi_{n,j,p}\in (D_{n,p})_r\cap {\mathcal F}$ for
any $j$. For a given $j$, this function is distinguished in
$(D_{n,p})_r$ by exactly $r$ (normalizing) conditions on matrices
$\xi$. Thus, $\psi_{n,j,p}$ is uniquely defined and the assertion
1 is proven.

Further on, observe that, for a given $n$, the space of all degree
$n$ functions in $\mathcal F$ is exactly a direct sum of its
subspaces $(D_{n,p})_r\cap {\mathcal F}$, $p=1,\ldots,N$. Hence,
its dimension equals to $rN$.
\end{proof}

\subsection{Fermion representations}\label{SS:ferm}${ }$

In this paragraph, we briefly outline the construction of a
\emph{fermion} $\D_\g$-module. It was  introduced in
\cite{rShf,ShMMJ} for the two-point situation ($N=1$). For the
multi-point situation considered here, the exposition is similar
up to a different definition of \KN\ basis (\refSS{knbb}).

Let $B$ be a rank $r$ degree $rg$ holomorphic vector bundle on
$\Sigma$, $\nabla$ a meromorphic connection on $B$ having simple
poles at most at the points of $A$ (hence, logarithmic --- see
\refS{orb}), $\tau$ an irreducible representation of $\g$ in the
finite-dimensional vector space $V_\tau$. Let $\Ga=\Ga(B)$ denote
the space of meromorphic sections of $B$ holomorphic except at
$P_1,\ldots,P_N,P_\infty$. Introduce $\Ga_{B,\tau}:=\Ga(B)\otimes
V_\tau$.

Define a $\Do_\g$-action on $\Ga_{B,\tau}$ as follows:

--- $\gb$-action: for any $x\in\g$, $A\in\A$, $s\in\Ga(B)$, $v\in
V_\tau$
\begin{equation}\label{E:gbact}
       (x\otimes A)(s\otimes v)=(A\cdot s)\otimes \tau(x)v ;
\end{equation}

--- $\L$-action: for any $e\in\L$, $s\in\Ga(B)$, $v\in V_\tau$
\begin{equation}\label{E:lac}
       e(s\otimes v)=\nabla_es\otimes v.
\end{equation}

\begin{proposition} Relations \refE{gbact}, \refE{lac} define a
Lie algebra representation of $\Do_\g$ in $\Ga_{B,\tau}$.
\end{proposition}
\begin{proof}
For the $\gb$ the claim is directly verified using \refE{gbact}.

By assumption, $\nabla$ is meromorphic, hence, flat. By flatness,
$\nabla_{[e,f]}=[\nabla_e,\nabla_f]$ for all $e,f\in\L$. Hence,
$\nabla$ defines a representation of $\L$ in $\Ga(B)$.

By definition of a connection, for any $s\in\Ga(B)$, $e\in\L$ and
$A\in\A$, we have $\nabla_e(As) = (e\ldot A)s+A\nabla_es$ where
$e\ldot A$ is the Lie derivative. Hence, $[\nabla_e,A]=e\ldot A$,
i.e.  the mapping $e+A\rightarrow\nabla_e+A$ gives rise to a
representation of $\D$ in $\Ga(B)$.
\end{proof}

Choose a \KN\ base $\{\psi_{n,p,j}\}$ in $\Ga(B)$ (see
\refSS{knbb}) and a weight base $\{v_a|1\le a\le \dim V_\tau\}$ in
$V_\tau$. Introduce $\psi_{n,p,j,a}=\psi_{n,p,j}\otimes v_a$.
Enumerate the elements $\psi_{n,p,j,a}$ linearly in the ascending
lexicographical order of the quadruples $(n,p,j,a)$. We write down
$\psi_M=\psi_{n,p,j,a}$ if $M=M(n,p,j,a)\in\Z$ is the number of
the quadruple $(n,p,j,a)$. Introduce the degree of $\psi_M$ by
$\deg \psi_M:=M$.

\begin{proposition}\label{L:agda}
With respect to the just introduced degree, $\Ga_{F,\tau}$ is an
almost-graded $\D_\g$-module.
\end{proposition}
For the two-point situation ($N=1$) the proof is given in
\cite{rShf,ShMMJ}). In case of general $N$ the proof is similar.

\medskip
Now we are in position to do the final step of the construction,
namely, to introduce the \emph{fermion space} corresponding to the
pair $(B,\tau)$ and define the $\Do_\g$-action in this space.

Consider the vector space $\Hwft$ generated over $\mathbb C$ by
the formal expressions ({\it semi-infinite wedge monomials}) of
the form $\Phi=\psi_{N_0}\wedge\psi_{N_1} \wedge\ldots$ where the
$\psi_{N_i}$ are the above introduced basis elements of
$\Ga_{F,\tau}$, the indices are strictly increasing, i.e.
$N_0<N_1<\ldots$, and $N_k=k+m$ for a suitable $m$ and all
sufficiently large $k$. The integer $m$ depends on the monomial
and following \cite{rKaRa} is called its {\it charge}. For a
monomial $\Phi$ of charge $m$, the degree of $\Phi$ is defined as
follows:
\begin{equation}\label{E:mdeg}
    \deg\Phi =\sum\limits_{k=0}^\infty (N_k-k-m).
\end{equation}
Observe that there is an arbitrariness in the enumeration of the
$\psi_{n,p,j,a}$'s for a fixed $n$; the just defined degree of a
monomial does not depend on this arbitrariness.

The monomials without requirement $N_0<N_1<\ldots$ are also
considered. It is assumed that they are antisymmetric with respect
to the order of their wedge co-multipliers.

We want to extend the action of $\Do_\g$ on $\Ga_{B,\tau}$ to an
action in  $\Hwft$. To this end, introduce the Lie algebra $\ainf$
of the infinite matrices with ``finitely many diagonals'' and
assign each element of $\D_\g$ with the matrix of its operator in
$\Ga(B)$ with respect to the basis $\psi_M$. Thus, we have
obtained the embedding of $\D_\g$ into $\ainf$. Then we use the
standard representation of $\ainf$ in $\Hwft$ due to V.Kac (for
example, see \cite{rKaRa}). The last step is absolutely standard
and was used in similar situations in \cite{rKNFb,rKNFc,rSDiss}.
We follow here the lines of \cite{rShf,ShMMJ}.

In turn, the construction due to V.Kac consists of two main steps,
namely applying the ``Leibnitz rule'' and the ``regularization''
if the result is not well-defined. The latter gives rise to a
representation of some central extension $\widehat{\Do_\g}$ of the
$\Do_\g$. Restrict the cocycle of this central extension to $\gb$
and denote so obtained cocycle by $\ga$.
\begin{proposition}\label{cokac} For $\g=\gl(l)$ and any
$x,y\in\g$, $f,g\in\A$
\[\ga(xf,yg)=\a\cdot\tr\,(xy)\ga^{(\A)}(f,g)
\]
(compare to \refE{eaff}, \refE{fung}) where $\a\in\C$.
\end{proposition}
The proof is the same as for Kac-Moody algebras \cite[p.97,
relation (9.16)]{rKaRa}.

Let $\Hwftm$ be the subspace of $\Hwft$ generated by the
semi-infinite monomials of charge $m$. These subspaces are
invariant under the action of $\widehat{\D_\g}$. Hence, for every
$m$ the space $\Hwftm$ itself is a $\widehat{\D_\g}$-module, and
$\Hwft=\bigoplus_{m\in\Z}\Hwftm$ as $\widehat{\D_g}$-module.

\begin{proposition}\label{P:fermalm}
Let $\Hwftm$ be the submodule of $\Hwft$ of charge $m$.
\newline
(a) With respect to the degree \refE{mdeg}, the homogeneous
subspaces $(\Hwftm)_k$ of degree $k$ are finite-dimensional. If
$k>0$ then $(\Hwftm)_k=0$.
\newline
(b) The cocycle $\gamma$ for  $\D_\g$  defined by the projective
representation is local.
\newline
(c) The $\Hwftm$ is an almost-graded $\widehat{\D_\g}$-module.
\end{proposition}
For $N=1$ the \refP{fermalm} is proven in \cite{rShf,ShMMJ}. For
the general case we refer to \cite{rSSpt2}.


\subsection{Sugawara representation}${ }$

Let $\gh$ be an affine algebra. A $\gh$-module $V$ is called
\emph{admissible} if any element of $V$ is annihilated by all
elements of $\gh_+$ of a sufficiently high degree.

Any admissible $\gh$-module can be canonically turned into an
$\Lh$-module by the \emph{Sugawara construction}. The
corresponding representation of $\Lh$ is called a \emph{Sugawara
representation}. For the conventional presentation of Sugawara
construction for Kac-Moody algebras, see \cite{rKaRa}. A
generalization onto the Krichever-Novikov algebras is given in
\cite{rKNFb,rBono,rSSS,rSSpt}.

Going over to the Sugawara construction, consider an admissible
$\gh$-module $V$ such that the central element $t$ operates as a
multiplication by a scalar $c$ (which is called a \emph{level} of
the representation).

For any $u\in\g$, $A\in \A$ denote by $\ u(A)\ $ the
representation operator of $u\otimes A$. Chose a basis $\ u_i,\
i=1,\ldots,\dim\g\ $ of $\g$ and the corresponding dual basis $\
u^i,\  i=1,\ldots,\dim\g$ with respect to the (fixed in advance)
invariant non-degenerate belinear form $(..|..)$.  We also denote
$u(A_{n,p})$ by $\ u(n,p)\ $ and $\sum_i u_i(n,p)u^i(m,q)$ by
$u(n,p)u(m,q)$ for short.

Define the higher genus {\it Sugawara operator} (also called {\it
Segal operator} or {\it energy-momentum tensor}) by
\begin{equation}\label{E:suga}
T(P):=\frac 12\sum_{n,m}\sum_{p,s}
\nord{u(n,p)u(m,s)}\w^{n,p}(P)\w^{m,s}(P)\ .
\end{equation}
where $\ \nord{....}\ $ denotes some normal ordering, $\w^{n,p}$
is the basis in the space of the 1-forms on $\Sigma$ introduced in
\refSS{knb}. Here, the summation indices $n,m$ run over $\Z$, and
$p,s$ over $\{1,\ldots,N\}$. The precise form of the normal
ordering is of no importance here.  For example, take the
following ``standard normal ordering'' ($x,y\in\g$)
\begin{equation}\label{E:normst}
\nord{x(n,p)y(m,r)}\ :=
\begin{cases} x(n,p)y(m,r),&n\le m
                                           \\
                            y(m,r)x(n,p),&n>m\
                      \end{cases}
\end{equation}
(for the discussion of normal orderings in case $g>0$ see
\cite{rKNFb,{rSSpt}}).

The expression $T(P)$ is considered as a formal series of
quadratic differentials on $Sigma$ with operator-valued
coefficients. Expanding it over the basis $\Omega^{k,r}$ of the
quadratic differentials (\refSS{knb}) we obtain
\begin{equation}\label{E:sugb}
  T(P)=\sum_k\sum_r L_{k,r}\cdot\Omega^{k,r}(P)\ ,
\end{equation}
with
\begin{equation}\label{E:sugc}
\begin{gathered}
  L_{k,r}=\cins T(P)e_{k,r}(P)=\frac 12\sum_{n,m}\sum_{p,s}
  \nord{u(n,p)u(m,s)}l_{(k,r)}^{(n,p)(m,s)},\\
  \text{where}\qquad
  \ l_{(k,r)}^{(n,p)(m,s)}:=\cins \w^{n,p}(P)\w^{m,s}(P)e_{k,r}(P)\ .
\end{gathered}
\end{equation}
A'priori, the operators $L_{k,r}$ are  infinite double sums. But
for given $k$ and $m$, the coefficient $l_{(k,r)}^{(n,p)(m,s)}$
will be non-zero only for finitely many $n$. This can be seen by
checking the residues of the elements appearing under the
integral. After applying the remaining infinite sum to a fixed
element $v\in V$, by the normal ordering and admissibility of the
representation only finitely many of the operators will operate
non-trivially on this element.

\begin{theorem}\label{T:suga} \cite{rSSS}
Let $\g$ be a finite dimensional either abelian or simple Lie
algebra and $2\k$ be the eigenvalue of its Casimir operator in the
adjoint representation. Let $V$ be  an admissible almost-graded
$\gh$-module of level $\ce$. If $\ce+\k\ne 0$ then the rescaled
modes
\begin{equation}\label{E:sugm}
L_{k,r}^*=\frac {-1}{2(\ce+\k)}\sum_{n,m}\sum_{p,s}
 \nord{u(n,p)u(m,s)}l_{(k,r)}^{(n,p)(m,s)}
\ ,
\end{equation}
of the Sugawara operator are well-defined operators on $V$ and
define an admissible representation of $\Lh$.
\end{theorem}
\begin{proposition}\label{P:sugalm}\cite{rSSpt}
The $V$ is an almost-graded $\Lh$-module under the Sugawara
action.
\end{proposition}
We call the $L_{k,r}^*$, resp. the $L_{k,r}$ the Sugawara
operators too. For $e=\sum_{n,p}a_{n,p}e_{n,p}\in\L$
($a_{n,p}\in{\mathbb C}$) we set $T[e]=\sum_{n,p}a_{n,p}L_{n,p}^*$
and obtain the representation $T$ of $\Lh$. It is called the {\it
Sugawara representation} of the Lie algebra $\L$ corresponding to
the given admissible representation $V$ of $\gh$.

By the Krichever-Novikov duality the Sugawara operator $T[e]$
assigned to the vector field
 $e\in\L$ can be given as
\begin{equation}\label{E:suge}
T[e]=\frac {-1}{\ce+\ka}\cdot \cins T(P)e(P).
\end{equation}

The following proposition expresses a fundamental property of the
Sugawara representation.
\begin{proposition}\label{P:fund}
For any reductive $\g$, $x\in\g$, $A\in\A$, $e\in\L$ we have
\begin{equation}\label{E:frel}
    [T[e],x(A)]=x(e\ldot A).
\end{equation}
\end{proposition}
\begin{proof} In case of a semi-simple or abelian $\g$, we refer
to \cite{rSSS,rSSpt} for a proof. Here, we give the proof for
$\g=\gl(l)$ which is only considered below. Actually, this is a
general proof in case of reductive $\g$.

In our case, $\g=\g_0\oplus\g_1$ where $\g_0$ is the center of
$\g$ consisting of diagonal matrices and $\g_1=\sln(l)$. Further
on, $\gb=\overline{\g_0}\oplus\overline{\g_1}$. Denote by
$x_0(A)$, $x_1(A)$ the restrictions of the representation $x(A)$
onto $\overline{\g_0}$, $\overline{\g_1}$, respectively. For
$x=x_0+x_1$, $x_0\in\g_0$, $x_1\in\g_1$  we can write
$x(A)=x_0(A)+x_1(A)$ without any conflict of notation. Moreover,
for any $A,B\in\A$, $x_0(A)$ and $x_1(B)$ commute because
$[x_0(A),x_1(B)]=[x_0,x_1]AB+\tr(x_0x_1)c\cdot id =0$.

Let $T_k$, $k=0,1$ be the Sugawara representation corresponding to
the representation $x_k$ of $\widehat{\g_k}$. Define
\[    T=T_0+T_1
\]
(in our case, this is equivalent to the definition \cite{KacB},
\cite[Rem. 10.2/10.3]{rKaRa}). Since operators of the
representations $T_0$ are expressed via $x_0(A)$'s and operators
of $T_1$ via $x_1(A)$'s, $T_0$ and $T_1$ commute, hence $T$ is a
representation of $\Lh$ and, moreover, for any $x_0\in\g_0$,
$x_1\in\g_1$, $A,B\in\A$, $e\in\L$
\[ [T_0[e],x_1(A)]=[T_1[e],x_0(B)]=0.
\]

For simple and abelian $\g$, the \refE{frel} is proven  in
\cite{rSSS} (see also \cite{rSSpt}). Hence,
\[ [T_0[e],x_0(A)]=x_0(e\ldot A),\quad [T_1[e],x_1(A)]=x_1(e\ldot A).
\]
Finally, we have
\[
 \begin{split}
    [T[e],x(A)]
    &=[T_0[e]+T_1[e], x_0(A)+x_1(A)] = [T_0[e], x_0(A)] +
      [T_1[e], x_1(A)] \\
    &= x_0(e.A)+x_1(e.A)= x(e.A).
  \end{split}
\]
\end{proof}


\section{Casimirs, semi-casimirs, Hitchin integrals}\label{S:cas}

Casimir operators (casimirs, laplacians) are the most important
invariant operators of Lie algebras. They are closely related to
all applications of Lie algebras and their representations
including integrable systems, the theory of special functions, and
many others. Description of casimirs is one of the central
problems of the representation theory. The second order casimirs
are of special interest in all these questions. In what follows,
``casimir" always means ``second order casimir".

Following the lines of \cite{ShMMJ}, we classify here the casimirs
for Krichever-Novikov algebras in case $N=1$ and $\g=\gl(r)$
($\g=\sln(r)$). We introduce the more general operators which we
call \emph{semi-casimirs}, investigate their relation to the
moduli space of Riemann surfaces and interpret them as a
quantization of the Hitchin intefrals. We delay the case $N>1$ for
the future publication \cite{rSSpt2}.

\subsection{Classification of casimirs}${ }$

Let $V$ be any admissible representation of $\Dh_g$. Its
restriction to $\gh$ is an admissible representation of the
latter.

Let $\ \eh\ $ denote the operator of representation of an
$e\in\L$, and $T[e]$ be the Sugawara operator of $e$ as introduced
above. Consider operators of the form $\Delta_e:=\eh -T[e]$.

In the following, we consider also the completed vector field
algebra $\Lc$ consisting of the infinite sums of the form
\begin{equation}\label{E:form}
    e=\sum\limits_{n= n_0}^\infty a_ne_n,\qquad
    a_n\in{\mathbb C}.\quad
\end{equation}
where, for $N=1$, $\ \{e_n=e_{n,1}\}$ is the Krichever-Novikov
base in $\L$; similarly $\{A_n=A_{n,1}\}$ is the Krichever-Novikov
base in $\A$. By admissibility, for any fixed vector $v$ both
$\eh\, v$ and $T[e]\,v$ are well-defined elements of $V$ even if
$e$ is of the form \refE{form}.

Let $\At_-\subset\A$ be the subspace spanned by all $A_k$, $k<0$.
\begin{definition}\label{D:semicas}
$ $\newline (a)   We call an operator of the form $\Delta_e$ a
{\it casimir} if $[\Delta_e,x(A)]=0$ for any $A\in\A$  and
$x\in\g$.
\newline (b) We call an operator of the form $\Delta_e$
a {\it semi-casimir} if $[\Delta_e,x(A)]=0$ for any $A\in\At_-$
and $x\in\g$.
\newline
\end{definition}
\begin{proposition}\label{P:eact}
Let $\g=\mathfrak{gl}(r)$ and $V$ be an admissible almost-graded
$\widehat{\D_\g}$-module such that the restriction of its cocycle
on $\gb$ is $\L$-invariant. Then for any $x\in\mathfrak{gl}(r)$
and $e\in\L$
\begin{equation}\label{E:comm2}
    [\Delta_e,x(A)]=\l(x)\ga^{(m)}(e,A)\cdot id,
\end{equation}
where $\ga^{(m)}$ is given by \refE{mixg} and
$\l(x)=r^{-1}\tr(x)$.
\end{proposition}
The \refP{eact} was first formulated in \cite{ShMMJ} under
different assumption on cocycles. The requirement of
$\L$-invariance (see \refSS{cocyc}) is proposed by
M.Schlichenmaier.
\begin{proof}
Since $V$ is an admissible almost-graded $\widehat{\D_\g}$-module,
its cocycle is local. By the classification of local
$\L$-invariant cocycles (\refSS{cocyc}) we have
\begin{equation}\label{E:comm1}
    [\eh,x(A)]=x(e\ldot A)+\l(x)\ga^{(m)}(e,A)\cdot id.
\end{equation}
Applying the \refP{fund} completes the proof.
\end{proof}

It follows from the \refP{eact} that for any $x\in\sln(r)$ (and
$A\in\A$) $\ [\Delta_e,xA]=0$. Thus, the cocycle $\ga^{(m)}$ is
the only obstacle for $\Delta_e$ to be a casimir; $\Delta(e)$ is a
casimir if and only if
\begin{equation}\label{E:sys}
\ga^{(m)}(A_k,e)=0, \quad \text{for any}\quad k\in\mathbb{Z}.
\end{equation}
Replacing here $e$ by its expression \refE{form}, obtain the
following linear system of equations on the coefficients
$\{a_n\}$:
\begin{equation}\label{E:syse}
\sum_{m\ge m_0}a_m\ga^{(m)}(A_{-k},e_m)=0, \quad \text{for
all}\quad k\in\mathbb{Z},\ k\ne 0.
\end{equation}
The further investigation of casimirs is based on the fact that
the system \refE{syse} is triangular and, for its diagonal
elements, in a generic situation, we have $\ga(A_{-k},e_k)\ne 0$,
for any $k\in\mathbb{Z}$, $\k\ne 0$. Since the equation for $k=0$
is missing, we obtain the 1-dimensional space of solutions
$\{a_n\}$ where $a_n=0$, $n<0$ and $a_n$ express via $a_0$, $n>0$.
This way, the description of casimirs can be completed in case of
the fermion modules. We formulate here only a final result and
refer to \cite{ShMMJ} for the proofs.

\begin{theorem}
For any fermion representation such that its cocycle $\ga$
satisfies the above condition of genericity, and any connection
$\nabla$ (involved via \refE{lac}), there exists exactly one (up
to a scalar factor) casimir. The corresponding vector field has a
simple zero at $P_+$.
\end{theorem}

\subsection{Semi-casimirs, coinvariants, moduli spaces }\label{SS:sscm}${ }$

Observe that for a vector field $e$ giving a semi-casimir one has
the system of linear equations similar to \refE{syse} but only for
$k>0$:
\begin{equation}\label{E:con2s}
\ga(A_{-k},e)=0,\quad \text{for any $k\in\mathbb{Z}$, $k>0$}.
\end{equation}
Thus, the coefficients $a_m$ with $m\le 0$ turn out to be
independent and all the others express via them. Let
$\Lt_-\subset\L$ be the subspace spanned by $\{e_k\colon k\le
0\}$.  Introduce the map $\Gamma\colon\Lt_-\mapsto\L$ as follows:
take $e\in\Lt_-$ and represent it in the form \refE{form}; then
substitute the corresponding $a_m$ ($m\le 0$) into \refE{con2s}
and calculate $a_m$, $m>0$. Denote by $\Gamma(e)$ the vector field
which corresponds to the full set of $a_m$'s.
\begin{lemma}\label{L:dss}\cite{ShMMJ}
The space of semi-casimirs coincides with $\Delta(\Gamma(\Lt_-))$.
It is spanned by the elements $\Delta(\Gamma(e_k))$, where $k\le
0$. It is isomorphic to $\L_-$ as a linear space.
\end{lemma}

Denote the subalgebra $\gh_{reg}$ introduced in \refSS{affine} by
$\g_r$, for short, and call it the \emph{regular subalgebra}. The
space of \emph{co-invariants} of $\g_r$ is defined as a quotient
space $V/U(\g_r)V$.

The semi-casimirs are defined in such way that they commute with
$U(\g_r)$, hence they are  well-defined on the space of
coinvariants of the subalgebra $\g_r$.

For $e\in\Lt_-$, let $\overline{\Delta}(e)$ be the operator
induced by $\Delta(\Gamma(e))$ on coinvariants.  The map
$\overline{\Delta}$ is defined on $\Lt_-$ and by \refL{dss} its
image is the space $C^s_2$ of semi-casimirs considered as
operators on the space of coinvariants.

Our next step is to show that only a finite number of basis
semi-casimirs are nonzero on coinvariants and, for a proper moduli
space of Riemann surfaces, establish the correspondence between
its tangent space and the space of semi-casimirs (considered on
coinvariants).

\begin{lemma}\label{L:kern}\cite{ShMMJ}
For a fermion representation $V$ there exists such
$p\in\mathbb{Z}_+$ that
$\L^{(p)}_-\subseteq\ker\overline{\Delta}$.
\end{lemma}
\begin{proof} It is proven in course of the proof of \cite[Lemma 4.9]{ShMMJ}
that such $p\in\Z$ exists that $\Delta(e)V\subseteq U(\g_r)V$ for
any $e\in\L^{(p)}_-$. Actually, we need the same for
$\Delta(\Gamma(e))$. Here, we only want to complete the proof of
\cite[Lemma 4.9]{ShMMJ} with this correction.

Let $V^{(q)}\subset V$ denote the subspace generated by all
elements of degree less or equal to $q$. In fact, it is proven in
course of the proof of \cite[Lemma 4.9]{ShMMJ} that for any
$q\in\Z$ such $p$ exists that $\Delta(e)V\subseteq V^{(q)}$ for
any $e\in\L^{(p)}_-$. Let $U_-$ be the subspace of the universal
enveloping algebra of $\gh$ generated by the basis elements of the
non-positive degree. Choose such $q$ that $U_-V^{(q)}\subseteq
U(\g_r)V$.

Take an arbitrary $e\in\L^{(p)}_-$. Observe that $\Gamma(e)=e+e_+$
where $e_+\in\L_+$. For any $v\in V$, we have $v=uv_0$ where $u\in
U_-$ and $v_0$ is the vacuum vector. Further on, $\Delta(e+e_+)$
is a semi-casimir, hence $\Delta(e+e_+)v=u\Delta(e+e_+)v_0$. By
\cite[Lemmas 3.2,\ 3.4]{ShMMJ}, $\Delta(e_+)v_0=0$, hence
$u\Delta(e+e_+)v_0=u\Delta(e)v_0$. Since $\Delta(e)v_0\subseteq
V^{(q)}$ and $uV^{(q)}\subseteq U(\g_r)V$, we have
$u\Delta(e)v_0\subseteq U(\g_r)V$, hence $\Delta(\Gamma(e))v\in
U(\g_r)V$.
\end{proof}

Let $\cM _{g,2}^{(p)}$ be the moduli space of curves of genus $g$
with two marked points $P_\pm$, fixed $1$-jet of local coordinate
at $P_+$ and fixed $p$-jet of local coordinate at $P_-$.  There is
a canonical mapping $\theta:\L\mapsto T_\Sigma\cM_{g,2}^{(p)}$
which goes back to \cite{Kon} and is based on the Kodaira-Spencer
theory. The cohomological and geometrical versions of this mapping
are given in \cite{rSSpt} and \cite{GrOr}, respectively (see
Introduction).

Let $\tilde\theta$ denote the restriction of $\theta$ onto the
subspace $\Lt_-$.  Let $V$ be a fermion representation of
$\widehat{\D_\g}$ and $\ga_V$ be its cocycle.  Let also
$C_2^s=C_2^s(V)$ denote the second order semi-casimirs of $\gh$ in
the representation $V$. We assume semi-casimirs to be restricted
onto coinvariants. The following theorem establishes a natural
mapping of the tangent space at $\Sigma\in\cM_{g,2}^{(p-1)}$ onto
the space of semi-casimirs in the coinvariants on $\Sigma$
considered as a punctured Riemann surface.

\begin{theorem}\label{T:mod}\cite{ShMMJ}
Take $p$ as in \refL{kern}.
\begin{enumerate}
\item[$1^\circ$.]  The mapping $\tilde\theta\colon \Lt_-\mapsto
T_\Sigma\cM_{g,2}^{(p-1)}$ is surjective and\/
$\ker\tilde\theta=\L^{(p)}_-$.

\item[$2^\circ$.]  For such $V$ that $\ga_V(A_{-k},e_k)\ne 0$ for any $k\in{\mathbb Z}_+$,
the mapping $\overline{\Delta}\colon \Lt_-\mapsto
C_2^s(V)$ is surjective and
$\L^{(p)}_-\subseteq\ker\overline{\Delta}$.

\item[$3^\circ$.] The mapping $\overline{\Delta}\circ\tilde\theta^{-1}\colon
T_\Sigma\cM_{g,2}^{(p-1)} \mapsto C_2^s(V)$ is well-defined and
surjective.
\end{enumerate}
\end{theorem}

\subsection{Quantization of the second order
              Hitchin integrals}

In this section, we show how the semi-casimirs appear in course of
operator quantization of the second order Hitchin integrals. We do
not give any mathematical setting the problem of quantization
here. We only show what happens if one follows some conventional
recipes.

Let $\phi$ is a Higgs field (mathematically, an arbitrary
$\g$-valued Krichever-Novikov 1-form on the Riemann surface in
question) and $\{\Omega^i\}$ is a base of the cotangent space to
${\mathcal M}_{g,2}^{(p-1)}$ realized as a certain space of
Krichever-Novikov quadratic differentials. We introduce the
\emph{second order Hitchin integrals} $\chi_i$'s by the expansion
$\tr\,\phi^2=\sum \chi_i\Omega^i$. Observe that this is only a
part of \emph{generalized} second order Hitchin integrals of
\refS{orb}, namely the part which contains no additional
functional factors.

As a first step of quantization, we replace $\phi$ by its
operator, thus we obtain the current $I$ which is an arbitrary
Krichever-Novikov 1-form (on the Riemann surface) having values in
the representation operators of $\gh$. Therefore, $I=\sum
u_k\omega^k$ where $\omega^k$ are the basis Krichever-Novikov
1-forms, $u_k$ are the operator-valued coefficients ($k\in{\Bbb
Z}$). Further on, the $\phi^2$ should be replaced by $\nord{I^2}$.
By definition of the Wess-Zumino-Witten-Novikov theory,
$\tr\nord{I^2}$ exactly equals to the energy-momentum tensor $T$
(introduced by \refE{suga}). The trace remains to be
"finite-dimensional", like in the classic situation, which means
that it is linear over the function algebra $\A$. The expansion
$T=\sum L_i\Omega^i$ (cf. \refE{sugb}) is the quantum analog of
the above expansion $\tr\,\phi^2=\sum \chi_i\Omega^i$. We will
consider the normalized form $-T(e_i)$ of an operator $L_i$ (see
Section~4.3). What is usually being done to compensate a normal
ordering, is adding certain cartanian elements to the normal
ordered quantity. For example, for Kac-Moody algebras the vector
field $e_0=z\frac{\partial}{\partial z}$ is considered and
$\widehat{e_0}$ is being added which leads to the casimir
$\widehat{z\frac{\partial}{\partial
z}}-T(z\frac{\partial}{\partial z})$. Applying this idea to an
arbitrary $e_i$, we come to the operators of the form
$\Delta_i=\widehat{e_i}-T(e_i)$. Since there is only a finite
number of (independent) Hitchin integrals and the infinite set of
$\Delta_i$'s, we formulate a selection rule for them. First, we
propose to consider certain linear combinations of $\Delta_i$'s,
namely those which are semi-casimirs. Thus, we replace $e_i$ by
$\Gamma(e_i)$ (where $\Gamma$ is introduced in \refSS{sscm}).
Since there is no canonical choice for the base $\{\Omega_i\}$,
the replacement $e_i$ by $\Gamma(e_i)$ can be achieved by
adjustment of this base. Second, we select only those
semi-casimirs which induce nontrivial operators on conformal
blocks. Thus, for each $i$ we consider $\Delta(\Gamma(e_i))$ as a
quantization of $\chi_i$. Then by Theorem~4.2 we obtain the
natural mapping of the $\chi_i$'s to the $\Delta(\Gamma(e_i))$'s.

Since the classic Hitchin integrals are in involution, the
corresponding quantum quantities must commute, at least, up to a
scalar. Let us show that this holds for our quantization of
Hitchin integrals.
\begin{proposition} In the space of a fermion representation, for any $e,f\in\L$,
\[  [\Delta(e),\Delta(f)]=\l(e,f)\cdot id
\]
where $\l$ is a bilinear form on  $\L$.
\end{proposition}
\begin{proof}
For any $u\in\g$, $A\in\A$ we have
\begin{align*}
    [[\Delta(e),\Delta(f)],u(A)]=&
    [[\Delta(e),u(A)],\Delta(f)]+[\Delta(e),[\Delta(f),u(A)] \\
    =&[(\tr\, u)\ga(e,A)\cdot id,\Delta(f)]+[\Delta(f),(\tr\,
    u)\ga(f,A)\cdot id]\\
    =&0.
\end{align*}
where $\gamma$ is a mixing cocycle. Thus, $[\Delta(e),\Delta(f)]$
is an automorphism of the $\gh$-module. For modules with a unique
highest vector (for example, the fermion modules) all
endomorphisms are scalar operators.
\end{proof}


\begin{thebibliography}{10}

\bibitem{Arn}
V.I. Arnold and Ju.S. Ilyashenko.
\newblock {\em Ordinary differential equations}, volume~1 of {\em Dynamical
  Systems 1. Sovremennye problemy matematiki. Fundamental'nye napravleniya.
  {D}.{V}.{A}nosov, {V}.{I}.{A}rnold, editors}.
\newblock M.: VINITI, 1985.

\bibitem{Bol}
A.A. Bolibrukh.
\newblock {\em 21-th problem of Hilbert for linear Fuchs systems}, volume 206
  of {\em Proceedings of Steklov Mathematical Institute}.
\newblock 1994.

\bibitem{rBono}
L.~Bonora, M.~Rinaldi, J.~Russo, and K.~Wu.
\newblock The {S}ugawara construction on genus $g$ {R}iemann surfaces.
\newblock {\em Phys. Lett. B}, 208:440--446, 1988.

\bibitem{Don}
S.~Donaldson.
\newblock Twisted harmonic maps and the self-duality equations.
\newblock {\em Proc. Lond. Math. Sci.}, (55)(1):127 -- 131, 1987.

\bibitem{GrOr}
P.G. Grinevich and A.Yu. Orlov.
\newblock Flag spaces in {KP} theory and {V}irasoro action on $det\ d$ and
  {S}egal-{W}ilson $\tau$-function.
\newblock {\em Preprint CLNS 945/89, Cornell University}, 1989.

\bibitem{rHsb}
N.~Hitchin.
\newblock Stable bundles and integrable systems.
\newblock {\em Duke Math. Journ.}

\bibitem{Had}
N.~Hitchin.
\newblock The self-duality equations on a riemann surface.
\newblock {\em Proc. Lond. Math. Sci.}, 55(1):59--126, 1987.

\bibitem{KacB}
V.G. Kac.
\newblock {\em Infinite dimensional {L}ie algebras}.
\newblock Cambridge Univ. Press, Cambridge, 1990.

\bibitem{rKaRa}
V.G. Kac and A.K. Raina.
\newblock {\em Highest Weight Representations of Infinite Dimensional Lie
  Algebras}, volume~2 of {\em Adv. Ser. in Math. Physics}.
\newblock World Scientific, 1987.

\bibitem{Kon}
M.~L. Kontsevich.
\newblock The {V}irasoro algebra and {T}eichm\"uller spaces.
\newblock {\em Funktsional. Anal. i Prilozhen. (Russian)}, 21(2):78--79, 1987.

\bibitem{rKNU}
I.M. Krichever and S.P. Novikov.
\newblock Holomorphic bundles on algebraic curves and nonlinear equations.
\newblock {\em Uspekhi Mat. Nauk}, (35)(6):47--68, 1980.

\bibitem{rKNFa}
I.M. Krichever and S.P. Novikov.
\newblock Algebras of {V}irasoro type, {R}iemann surfaces and structures of the
  theory of solitons.
\newblock {\em Funktional Anal. i. Prilozhen.}, 21(2):46--63, 1987.

\bibitem{rKNFb}
I.M. Krichever and S.P. Novikov.
\newblock {V}irasoro type algebras, {R}iemann surfaces and strings in
  {M}inkowski space.
\newblock {\em Funktional Anal. i. Prilozhen.}, 21(4):47--61, 1987.

\bibitem{rKNFc}
I.M. Krichever and S.P. Novikov.
\newblock Algebras of {V}irasoro type, energy-momentum tensors and
  decompositions of operators on {R}iemann surfaces.
\newblock {\em Funktional Anal. i. Prilozhen.}, 23(1):46--63, 1989.

\bibitem{rKNR2p}
I.M. Krichever and S.P. Novikov.
\newblock Holomorphic bundles and commuting difference operators. {T}wo-point
  constructions.
\newblock {\em Usp. Mat. Nauk}, 55(3):181--182, 2000.

\bibitem{rRDS}
A.~Ruffing, Th. Deck, and M.~Schlichenmaier.
\newblock String branchings on complex tori and algebraic representations of
  generalized {K}richever-{N}ovikov algebras.
\newblock {\em Lett. Math. Phys.}, 26:23--32, 1992.

\bibitem{rSad}
V.A. Sadov.
\newblock Bases on multipunctured {R}iemann surfaces and interacting strings
  amplitudes.
\newblock {\em Commun. Math. Phys.}, 136:585--597, 1991.

\bibitem{SchlMMJ}
M.~Schlichenmaier.
\newblock {H}igher genus affine algebras of {K}richever-{N}ovikov type.
\newblock {\em Moscow Math. Journ.}
\newblock To appear. Math.QA/0210360.

\bibitem{rSLc}
M.~Schlichenmaier.
\newblock Central extensions and semi-infinite wedge representations of
  {K}richever-{N}ovikov algebras for more than two points.
\newblock {\em Letters in Mathematical Physics}, 20:33--46, 1990.

\bibitem{rSLa}
M.~Schlichenmaier.
\newblock Krichever-novikov algebras for more than two points.
\newblock {\em Letters in Mathematical Physics}, 19:151--165, 1990.

\bibitem{rSLb}
M.~Schlichenmaier.
\newblock {K}richever-{N}ovikov algebras for more than two points: explicit
  generators.
\newblock {\em Letters in Mathematical Physics}, 19:327--336, 1990.

\bibitem{rSDiss}
M.~Schlichenmaier.
\newblock {\em Verallgemeinerte {K}richever - {N}ovikov {A}lgebren und deren
  {D}arstellungen}.
\newblock PhD thesis, Universit{\"{a}}t Mannheim, 1990.

\bibitem{rSDeg}
M.~Schlichenmaier.
\newblock Degenerations of generalized {K}richever-{N}ovikov algebras on tori.
\newblock {\em Journal of Mathematical Physics}, 34:3809--3824, 1993.

\bibitem{SchlCo}
M.~Schlichenmaier.
\newblock Local cocycles and central extensions for multi-point algebras of
  {K}richever-{N}ovikov type.
\newblock {\em J. Reine und Angewandte Mathematik}, 2001.
\newblock Math/0112116.

\bibitem{rSSS}
M.~Schlichenmaier and O.~Sheinman.
\newblock {S}ugawara construction and {C}asimir operators for
  {K}richever-{N}ovikov algebras.
\newblock {\em Jour. of Math. Science}, 92:3807--3834, 1998.

\bibitem{rSSpt2}
M.~Schlichenmaier and O.K. Sheinman.
\newblock {W}ess-{Z}umino-{W}itten-{N}ovikov theory, {K}nizhnik-{Z}amolodchikov
  equations, and {K}richever-{N}ovikov algebras, {II}.
\newblock In progress.

\bibitem{rSSpt}
M.~Schlichenmaier and O.K. Sheinman.
\newblock {W}ess-{Z}umino-{W}itten-{N}ovikov theory, {K}nizhnik-{Z}amolodchikov
  equations, and {K}richever-{N}ovikov algebras, {I}.
\newblock {\em Russian Math. Surv. (Uspekhi Math. Nauk).}, 54:213--250, 1999.
\newblock Math.QA/9812083.

\bibitem{rShea}
O.K. Sheinman.
\newblock {E}lliptic affine {L}ie algebras.
\newblock {\em Funktional Anal. i. Prilozhen.}, 24(3):210--219, 1992.

\bibitem{rSha}
O.K. Sheinman.
\newblock Affine {L}ie algebras on {R}iemann surfaces.
\newblock {\em Funktional Anal. i. Prilozhen.}, 27(4):54--62, 1993.

\bibitem{rShf}
O.K. Sheinman.
\newblock The fermion model of representations of affine {K}richever-{N}ovikov
  algebras.
\newblock {\em Funktional Anal. i. Prilozhen.}, 35(3), 2001.
\newblock Math.RT/0204178.

\bibitem{ShSD}
O.K. Sheinman.
\newblock {K}richever-{N}ovikov algebras and the self-duality equations on
  {R}iemann surfaces.
\newblock {\em Russian Math. Surv.}, 56(1):185--186, 2001.

\bibitem{ShMMJ}
O.K. Sheinman.
\newblock Second order casimirs for the affine {K}richever-{N}ovikov algebras
  $\widehat{\mathfrak{gl}}_{g,2}$ and $\widehat{\mathfrak{sl}}_{g,2}$.
\newblock {\em Moscow Math.J.}, 1(4):605--628, 2001.
\newblock Math.RT/0109001.

\bibitem{ShUMN}
O.K. Sheinman.
\newblock {S}econd order casimirs for the affine {K}richever-{N}ovikov algebras
  $\widehat{\frak{gl}}_{g,2}$ and $\widehat{\frak{sl}}_{g,2}$.
\newblock {\em Russian Math. Surv.}, 56:189--190, 2001.

\bibitem{rTur}
A.N. Tjurin.
\newblock Classification of vector bundles on an algebraic curve of arbitrary
  genus.
\newblock {\em Izvestia AN SSSR}, 29:657--688, 1965.

\end{thebibliography}

\end{document}